\documentclass{amsart}%
\usepackage{amsfonts}
\usepackage{graphicx}
\usepackage{amscd}%
\usepackage{amsmath}%
\usepackage{amssymb}
\newtheorem{theorem}{Theorem}
\theoremstyle{plain}

\newtheorem{corollary}{Corollary}

\newtheorem{lemma}{Lemma}

\newtheorem{proposition}{Proposition}

\numberwithin{equation}{section}
\numberwithin{theorem}{section}
\numberwithin{algorithm}{section}
\numberwithin{axiom}{section}
\numberwithin{case}{section}
\numberwithin{claim}{section}
\numberwithin{conclusion}{section}
\numberwithin{condition}{section}
\numberwithin{conjecture}{section}
\numberwithin{corollary}{section}
\numberwithin{criterion}{section}
\numberwithin{definition}{section}
\numberwithin{example}{section}
\numberwithin{exercise}{section}
\numberwithin{lemma}{section}
\numberwithin{notation}{section}
\numberwithin{problem}{section}
\numberwithin{proposition}{section}
\numberwithin{remark}{section}
\numberwithin{solution}{section}

\begin{document}
\title[Conformal covariant operators on the sphere]{On the higher order conformal covariant operators on the sphere }
\author{Fengbo Hang}
\address{Department of Mathematics, Michigan State University, East Lansing, MI 48824}
\email{fhang@math.msu.edu}

\begin{abstract}
We will show that in the conformal class of the standard metric $g_{S^{n}}$ on
$S^{n}$, the scaling invariant functional $\left(  \mu_{g}\left(
S^{n}\right)  \right)  ^{\frac{2m-n}{n}}\int_{S^{n}}Q_{2m,g}d\mu_{g}$
maximizes at $g_{S^{n}}$ when $n$ is odd and $m=\frac{n+1}{2}$ or $\frac
{n+3}{2}$. For $n$ odd and $m\geq\frac{n+5}{2}$, $g_{S^{n}}$ is not stable and
the functional has no local maximizer. Here $Q_{2m,g}$ is the $2m$th order $Q
$-curvature.

\end{abstract}
\maketitle

\section{Introduction}

Let $A_{g}$ be a differential operator on a $n$-dimensional Riemannian
manifold $\left(  M,g\right)  $. Recall that we say $A$ is conformally
covariant of bidegree $\left(  a,b\right)  $ if for any $u,w\in C^{\infty
}\left(  M\right)  $,
\[
A_{e^{2w}g}u=e^{-bw}A_{g}\left(  e^{aw}u\right)  .
\]
The most well known conformal covariant operators are the Laplacian operator
on a surface, which is of bidegree $\left(  0,2\right)  $ and the conformal
Laplacian operator on a manifold of dimension $n\geq3$,
\[
L_{g}=-\Delta_{g}+\frac{n-2}{4\left(  n-1\right)  }R,
\]
which is of bidegree $\left(  \frac{n-2}{2},\frac{n+2}{2}\right)  $. Here $R$
denotes the scalar curvature. They play important role in the study of
Gaussian curvature and scalar curvature. Besides these two examples, on a four
dimensional manifold, the fourth order Paneitz operator discovered in
\cite{P}, which is of bidegree $\left(  0,4\right)  $, has demonstrated its
importance in conformal geometry recently (cf. \cite{CGY}). According to
\cite{FG1,GJMS,Br}, there exists a sequence of conformal covariant operators
which contain the three examples above (see also \cite{FG2,GZ}). Indeed, if
$m$ is a positive integer such that either $n$ is odd and $n\geq3$, or $n$ is
even and $2m\leq n$, there exists a conformally covariant operator $P_{2m} $
of bidegree $\left(  \frac{n-2m}{2},\frac{n+2m}{2}\right)  $. Moreover, the
leading term of $P_{2m}$ is equal to $\left(  -\Delta_{g}\right)  ^{m}$ and on
$\mathbb{R}^{n}$ with standard metric, $P_{2m}=\left(  -\Delta\right)  ^{m} $.
It is an interesting fact pointed out by \cite{G,GH} that the condition
$2m\leq n$ is necessary when $n\geq4$ is even. For recent developments related
to these operators, one should refer to \cite{A,C} and the references therein.

In general, it seems very hard to have an explicit formula for the operators
$P_{2m}$. However, on $S^{n}$ with standard metric, $P_{2m}$ has a nice
expression as (see part (f) of theorem 2.8 in \cite{Br})%
\begin{equation}
P_{2m}=%
{\displaystyle\prod\limits_{i=0}^{m-1}}
\left(  -\Delta_{S^{n}}-\left(  i+\frac{n}{2}\right)  \left(  i-\frac{n}%
{2}+1\right)  \right)  .\label{eqP2m}%
\end{equation}
We note that one does not require $2m\leq n$ when $n$ is even in this special
case (see Section \ref{secP2mExpression} for more information).

Assume $n\geq3$, let $g_{S^{n}}$ be the standard metric on $S^{n}$, then we
know (see chapter V of \cite{SY})%
\begin{align*}
Y\left(  S^{n}\right)   & =\inf\left\{  \frac{\int_{S^{n}}R_{g}d\mu_{g}}%
{\mu_{g}\left(  S^{n}\right)  ^{\frac{n-2}{n}}}:g=\rho^{2}g_{S^{n}},\rho\in
C^{\infty}\left(  S^{n},\mathbb{R}\right)  ,\rho>0\right\} \\
& =n\left(  n-1\right)  \left(  \mu_{S^{n}}\left(  S^{n}\right)  \right)
^{\frac{2}{n}}\text{.}%
\end{align*}
That is, the functional minimizes at the standard metric. Here $R_{g}$ is the
scalar curvature of $g$ and $\mu_{g}$ is the measure associated with $g$.
Moreover $g$ is a critical metric if and only if $g=c\cdot\phi^{\ast}g_{S^{n}%
}$ for some positive number $c$ and Mobius transformation $\phi$, and all of
them are minimizers. In terms of the conformal Laplacian operator, we have%
\[
Y\left(  S^{n}\right)  =\inf\left\{  \frac{4\left(  n-1\right)  }{n-2}%
\frac{\int_{S^{n}}L_{S^{n}}u\cdot ud\mu_{S^{n}}}{\left(  \int_{S^{n}}%
u^{\frac{2n}{n-2}}d\mu_{S^{n}}\right)  ^{\frac{n-2}{n}}}:u\in C^{\infty
}\left(  S^{n},\mathbb{R}\right)  ,u>0\right\}
\]
and the minimizing value is reached at $u=1$, moreover, $u$ is a critical
point if and only if $u=c\cdot J_{\phi}^{\frac{n-2}{2n}}$ for some positive
number $c$ and Mobius transformation $\phi$. Here $J_{\phi}$ denotes the
Jacobian of $\phi$. Using the stereographic projection from $S^{n}%
\backslash\left\{  N\right\}  $ ($N$ is the north pole of $S^{n}$) to
$\mathbb{R}^{n}$ and the simple fact that for any $u\in C^{\infty}\left(
S^{n}\right)  $, $u\geq0$, we may find a sequence $u_{i}\in C_{c}^{\infty
}\left(  S^{n}\backslash\left\{  N\right\}  \right)  $, $u_{i}\geq0$ such that
$u_{i}\rightarrow u$ in $H^{1}\left(  S^{n}\right)  $, we see%
\begin{align*}
Y\left(  S^{n}\right)   & =\inf\left\{  \frac{4\left(  n-1\right)  }{n-2}%
\frac{\int_{\mathbb{R}^{n}}\left\vert \nabla\varphi\right\vert ^{2}%
d\mu_{\mathbb{R}^{n}}}{\left(  \int_{\mathbb{R}^{n}}\varphi^{\frac{2n}{n-2}%
}d\mu_{\mathbb{R}^{n}}\right)  ^{\frac{n-2}{n}}}:\varphi\in C_{c}^{\infty
}\left(  \mathbb{R}^{n}\right)  ,\varphi\geq0,\varphi\equiv\hspace
{-0.13in}\backslash\ 0\right\} \\
& =\inf\left\{  \frac{4\left(  n-1\right)  }{n-2}\frac{\int_{\mathbb{R}^{n}%
}\left\vert \nabla\varphi\right\vert ^{2}d\mu_{\mathbb{R}^{n}}}{\left(
\int_{\mathbb{R}^{n}}\varphi^{\frac{2n}{n-2}}d\mu_{\mathbb{R}^{n}}\right)
^{\frac{n-2}{n}}}:\varphi\in C_{c}^{\infty}\left(  \mathbb{R}^{n}\right)
,\varphi\equiv\hspace{-0.13in}\backslash\ 0\right\}  .
\end{align*}
In particular, the statement that the standard metric is a minimizer is
equivalent to the sharp Sobolev inequalities studied earlier in \cite{Au,T}.

When $n=2$, the parallel statement is that the standard metric $g_{S^{2}}$ has
the maximal determinant among all smooth metrics $g=e^{2u}g_{S^{2}}$ with
$\mu_{g}\left(  S^{2}\right)  =4\pi$ (see theorem 1 of \cite{OPS}). More
precisely, we have the Polyakov formula%
\[
\log\frac{\det^{\prime}\Delta_{e^{2u}g_{S^{2}}}}{\det^{\prime}\Delta_{S^{2}}%
}=-\frac{1}{12\pi}\int_{S^{2}}\left\vert \nabla u\right\vert ^{2}d\mu_{S^{2}%
}-\frac{1}{6\pi}\int_{S^{2}}ud\mu_{S^{2}}%
\]
for $u\in C^{\infty}\left(  S^{2}\right)  $ with $\int_{S^{2}}e^{2u}%
d\mu_{S^{2}}=4\pi$, and the Onofri inequality%
\[
\log\left(  \frac{1}{4\pi}\int_{S^{2}}e^{2u}d\mu_{S^{2}}\right)  \leq\frac
{1}{4\pi}\int_{S^{2}}\left\vert \nabla u\right\vert ^{2}d\mu_{S^{2}}+\frac
{1}{2\pi}\int_{S^{2}}ud\mu_{S^{2}}%
\]
for $u\in C^{\infty}\left(  S^{2}\right)  $.

When $2m<n$, the $Q$-curvature $Q_{2m,g}$ is given by%
\[
Q_{2m,g}=\frac{2}{n-2m}P_{2m,g}1,
\]
(cf. theorem 1.1 of \cite{Br}). We have%
\begin{align*}
& Y_{2m}\left(  S^{n}\right) \\
& =\inf\left\{  \frac{\int_{S^{n}}Q_{2m,g}d\mu_{g}}{\left(  \mu_{g}\left(
S^{n}\right)  \right)  ^{\frac{n-2m}{n}}}:g=\rho^{2}g_{S^{n}},\rho\in
C^{\infty}\left(  S^{n},\mathbb{R}\right)  ,\rho>0\right\} \\
& =\inf\left\{  \frac{2}{n-2m}\frac{\int_{S^{n}}P_{2m}u\cdot ud\mu_{S^{n}}%
}{\left(  \int_{S^{n}}u^{\frac{2n}{n-2m}}d\mu_{S^{n}}\right)  ^{\frac{n-2m}%
{n}}}:u\in C^{\infty}\left(  S^{n},\mathbb{R}\right)  ,u>0\right\} \\
& =\inf\left\{  \frac{2}{n-2m}\frac{\int_{\mathbb{R}^{n}}\left\vert
D^{m}\varphi\right\vert ^{2}dx}{\left(  \int_{\mathbb{R}^{n}}\varphi
^{\frac{2n}{n-2m}}dx\right)  ^{\frac{n-2m}{n}}}:\varphi\in C_{c}^{\infty
}\left(  \mathbb{R}^{n},\mathbb{R}\right)  ,\varphi\geq0,\varphi\equiv
\hspace{-0.13in}\backslash\ 0\right\}  .
\end{align*}
It follows from \cite{Lie,Lin,Lio,S,WX} that the standard metric $g_{S^{n}}$
is a minimizer (i.e. $u=1$ is a minimizer). The case $2m=n$ was treated in
\cite{B,CY,Lie}. When $2m>n$, for $Q_{2m,g}=\frac{2}{n-2m}P_{2m,g}1$ (assume
$n$ is odd, see theorem 1.1 in \cite{Br}), we have%
\begin{align*}
& Y_{2m}\left(  S^{n}\right) \\
& =\sup\left\{  \left(  \mu_{g}\left(  S^{n}\right)  \right)  ^{\frac{2m-n}%
{n}}\int_{S^{n}}Q_{2m,g}d\mu_{g}:g=\rho^{2}g_{S^{n}},\rho\in C^{\infty}\left(
S^{n},\mathbb{R}\right)  ,\rho>0\right\} \\
& =\frac{2}{n-2m}\inf\left\{  \left\vert u^{-1}\right\vert _{L^{\frac
{2n}{2m-n}}\left(  S^{n}\right)  }^{2}\int_{S^{n}}P_{2m}u\cdot ud\mu_{S^{n}%
}:u\in C^{\infty}\left(  S^{n},\mathbb{R}\right)  ,u>0\right\}  .
\end{align*}
We are motivated to ask the following question: Is the standard metric
$g_{S^{n}}$ a maximizer? Or equivalently: Is%

\begin{align}
& B_{2m}\left(  S^{n}\right) \label{eqmain}\\
& =\inf\left\{  \left\vert u^{-1}\right\vert _{L^{\frac{2n}{2m-n}}\left(
S^{n}\right)  }^{2}\int_{S^{n}}P_{2m}u\cdot ud\mu_{S^{n}}:u\in C^{\infty
}\left(  S^{n},\mathbb{R}\right)  ,u>0\right\} \nonumber
\end{align}
achieved at $u=1$? Due to the existence of negative power, this variational
problem is analytically different from the case $2m<n$. The answer to the
above question would shed some light on the understanding of the now still
mysterious $Q$-curvature with $2m>n$. We remark that the operator $P_{2}$ on
$S^{1}$ appears naturally in the study of self-similar solutions for the
anisotropic affine curve shortening problem (cf. \cite{ACW}). In particular,
the variational problem (\ref{eqmain}) indeed has $u=1$ as a minimizer in this
case (see proposition 1.3 in \cite{ACW}). Another interesting case of $P_{4}$
on $S^{3}$ was solved affirmatively in \cite{YZ} (see section 7 of \cite{HY}
for a different proof). The main aim of this note is to resolve all the
remaining cases. To state the result, we introduce some notations.

For $u,v\in C^{\infty}\left(  S^{n},\mathbb{R}\right)  $, we denote%
\begin{equation}
E_{2m}\left(  u,v\right)  =\int_{S^{n}}P_{2m}u\cdot vd\mu_{S^{n}%
}.\label{eqE2m1}%
\end{equation}
By integration by parts and the standard approximation argument, we know
$E_{2m}$ has a unique bounded symmetric bilinear extension to $H^{m}\left(
S^{n}\right)  \times H^{m}\left(  S^{n}\right)  $. Here $H^{m}\left(
S^{n}\right)  =W^{m,2}\left(  S^{n}\right)  $. Denote%
\begin{equation}
E_{2m}\left(  u\right)  =E_{2m}\left(  u,u\right)  \quad\text{for }u\in
H^{m}\left(  S^{n}\right)  .\label{eqE2m2}%
\end{equation}
Assume $2m>n$, then $H^{m}\left(  S^{n}\right)  \subset C\left(  S^{n}\right)
$, hence we may define%
\begin{equation}
V_{m}=\left\{  u\in H^{m}\left(  S^{n}\right)  :u>0\text{ everywhere}\right\}
.\label{eqVm}%
\end{equation}
Denote%
\begin{equation}
I_{2m}\left(  u\right)  =\left\vert u^{-1}\right\vert _{L^{\frac{2n}{2m-n}%
}\left(  S^{n}\right)  }^{2}E_{2m}\left(  u\right)  ,\quad B_{2m}\left(
S^{n}\right)  =\inf_{u\in V_{m}}I_{2m}\left(  u\right)  .\label{eqI2m}%
\end{equation}

\begin{theorem}
\label{thmmain}Let $P_{2m}$ be the $2m$th order conformal covariant operator
on $S^{n}$ (see (\ref{eqP2m})), $E_{2m}$, $V_{m}$, $I_{2m}$ and $B_{2m}\left(
S^{n}\right)  $ be defined as in (\ref{eqE2m1}), (\ref{eqE2m2}), (\ref{eqVm})
and (\ref{eqI2m}).

\begin{itemize}
\item[(i)] If $n$ is odd and $m=\frac{n+1}{2}$, then for any $u\in
H^{\frac{n+1}{2}}\left(  S^{n}\right)  $, $u>0$,%
\begin{align*}
\left\vert u^{-1}\right\vert _{L^{2n}\left(  S^{n}\right)  }^{2}E_{n+1}\left(
u\right)   & \geq E_{n+1}\left(  1\right)  \left(  \mu_{S^{n}}\left(
S^{n}\right)  \right)  ^{1/n}\\
& =-\frac{\left(  2n\right)  !}{2^{2n+1}\cdot n!}\left(  \frac{2\pi
^{\frac{n+1}{2}}}{\left(  \frac{n-1}{2}\right)  !}\right)  ^{\frac{n+1}{n}}.
\end{align*}
Moreover, all the minimizers of $I_{n+1}$ on $V_{\frac{n+1}{2}}$ are of the
form $cJ_{\phi}^{-\frac{1}{2n}}$ for some $c>0$ and Mobius transformation
$\phi$. Here $J_{\phi}$ is the Jacobian of $\phi$.

\item[(ii)] If $n$ is odd and $m=\frac{n+3}{2}$, then for any $u\in
H^{\frac{n+3}{2}}\left(  S^{n}\right)  $, $u>0$,%
\begin{align*}
\left\vert u^{-1}\right\vert _{L^{\frac{2n}{3}}\left(  S^{n}\right)  }%
^{2}E_{n+3}\left(  u\right)   & \geq E_{n+3}\left(  1\right)  \left(
\mu_{S^{n}}\left(  S^{n}\right)  \right)  ^{3/n}\\
& =\frac{3\cdot\left(  2n+1\right)  !}{2^{2n+3}\cdot n!}\left(  \frac
{2\pi^{\frac{n+1}{2}}}{\left(  \frac{n-1}{2}\right)  !}\right)  ^{\frac
{n+3}{n}}.
\end{align*}
Again, all the minimizers of $I_{n+3}$ in $V_{\frac{n+3}{2}}$ are of the form
$cJ_{\phi}^{-\frac{3}{2n}}$ for some $c>0$ and Mobius transformation $\phi$.

\item[(iii)] If $n$ is odd and $m\geq\frac{n+5}{2}$, then $I_{2m}$ has no
local minimizer in $V_{2m}$. Indeed, all the critical points of $I_{2m}$ in
$V_{m}$ are unstable.

\item[(iv)] If $n$ is even and $m>\frac{n}{2}$, then $P_{2m}\geq0$, moreover%
\[
\ker P_{2m}=\left\{  \left.  p\right\vert _{S^{n}}:p\text{ is a polynomial on
}\mathbb{R}^{n+1}\text{ with }\deg p\leq m-\frac{n}{2}\right\}  .
\]
In particular, $I_{2m}$ minimizes at $u=1$ and all the minimizers of $I_{2m}$
in $V_{m}$ are of the form $\left.  p\right\vert _{S^{n}}$, where $p$ is a
polynomial of degree less than or equal to $m-\frac{n}{2}$ and $\left.
p\right\vert _{S^{n}}>0$.
\end{itemize}
\end{theorem}

To get a feeling of the inequalities proved here, note that for $P_{2}$ on
$S^{1}$, what we have got is%
\begin{equation}
\int_{S^{1}}u^{-2}d\theta\int_{S^{1}}\left(  u_{\theta}^{2}-\frac{1}{4}%
u^{2}\right)  d\theta\geq-\pi^{2}\label{eqP2S1}%
\end{equation}
for $u\in C^{\infty}\left(  S^{1}\right)  $, $u>0$. This inequality was proved
earlier in proposition 1.3 of \cite{ACW}.

For $P_{4}$ on $S^{1}$, what we have got is
\begin{equation}
\left(  \int_{S^{1}}u^{-2/3}d\theta\right)  ^{3}\int_{S^{1}}\left(
u_{\theta\theta}^{2}-\frac{5}{2}u_{\theta}^{2}+\frac{9}{16}u^{2}\right)
d\theta\geq9\pi^{4}\label{eqP4S1}%
\end{equation}
for all $u\in C^{\infty}\left(  S^{1}\right)  $, $u>0$. It is interesting to
note that if we take $u=\sin\theta$, then the left hand side of (\ref{eqP4S1})
is a finite negative number. Hence the condition $u>0$ is crucial for the
validity of (\ref{eqP4S1}).

The article will be written as follows: In Section \ref{sec2}, we will give an
elementary argument for the expression of $P_{2m}$ on $S^{n}$ (cf.
(\ref{eqP2m})) and its invariant property under the Mobius transformation,
then we will discuss when a Sobolev function can be approximated by functions
vanishing near a given point. After these preparations, we shall prove Theorem
\ref{thmmain} in Section \ref{sec3}. In Section \ref{sec4} we give a somewhat
different argument for Theorem \ref{thmmain} based on the barycenter analysis.
In the last section we make some remarks concerning the proof of
(\ref{eqP4S1}) given in the recent preprint \cite{NZ}.

\textbf{Acknowledgment}: The research of the author is supported by National
Science Foundation Grant DMS-0209504. We would like to thank Paul Yang for
valuable discussions.

\section{Some preparations\label{sec2}}

First let us fix the notation for stereographic projection from the punctured
sphere to the Euclidean space which we will use later. For any $\xi\in S^{n}$,
let
\[
\xi^{\perp}=\left\{  z\in\mathbb{R}^{n+1}:z\cdot\xi=0\right\}  .
\]
For every $z\in\mathbb{R}^{n+1}$,
\[
z=z^{\prime}+t\xi=\left(  z^{\prime},t\xi\right)  ,\quad z^{\prime}\in
\xi^{\perp},t\in\mathbb{R}.
\]
The stereographic projection is
\[
\pi_{\xi}:S^{n}\backslash\left\{  \xi\right\}  \rightarrow\xi^{\perp
}:z=\left(  z^{\prime},t\xi\right)  \mapsto\frac{z^{\prime}}{1-t},
\]
its inverse is
\[
\pi_{\xi}^{-1}:\xi^{\perp}\rightarrow S^{n}\backslash\left\{  \xi\right\}
:x\mapsto\left(  \frac{2x}{\left\vert x\right\vert ^{2}+1},\frac{\left\vert
x\right\vert ^{2}-1}{\left\vert x\right\vert ^{2}+1}\xi\right)  .
\]
We have%
\[
\left(  \pi_{\xi}^{-1}\right)  ^{\ast}g_{S^{n}}=\frac{4}{\left(  1+\left\vert
x\right\vert ^{2}\right)  ^{2}}\sum_{i=1}^{n}dx_{i}\otimes dx_{i},
\]
here $x_{1},\cdots,x_{n}$ is the coordinate on $\xi^{\perp}$ with respect to
any fixed orthonormal frame of $\xi^{\perp}$.

For $\lambda>0$, we have a Mobius transformation $\sigma_{\xi,\lambda}\left(
\zeta\right)  =\pi_{\xi}^{-1}\left(  \lambda\pi_{\xi}\left(  \zeta\right)
\right)  $, it satisfies%
\[
\sigma_{\xi,\lambda}^{\ast}g_{S^{n}}=\frac{\lambda^{2}\left(  1+\left\vert
\pi_{\xi}\right\vert ^{2}\right)  ^{2}}{\left(  1+\lambda^{2}\left\vert
\pi_{\xi}\right\vert ^{2}\right)  ^{2}}g_{S^{n}}.
\]

\subsection{An elementary argument to derive the expression of $P_{2m}$ and
its properties\label{secP2mExpression}}

In this subsection, we will derive the expression of $P_{2m}$ on $S^{n}$ by an
elementary induction argument. Along the way, we shall also derive the
transformation law of $P_{2m}$ which we will use later. In principle, it makes
the proof of inequalities in Theorem \ref{thmmain} self-contained. We should
point out that the expression of $P_{2m}$ on $S^{n}$ was explicitly written
down in the part (f) of theorem 2.8 in \cite{Br}.

\begin{lemma}
\label{lemidentity1}Let $u$ be a smooth function on a domain in $\mathbb{R}%
^{n}$, and $m$ be a nonnegative integer, then%
\begin{align*}
& \Delta\left[  \left(  \frac{1+\left\vert x\right\vert ^{2}}{2}\right)
^{m+1}\Delta^{m}u\right]  +m\left(  m+1\right)  \left(  \frac{1+\left\vert
x\right\vert ^{2}}{2}\right)  ^{m-1}\Delta^{m}u\\
& =\left(  \frac{1+\left\vert x\right\vert ^{2}}{2}\right)  ^{m}\Delta
^{m+1}\left(  \frac{1+\left\vert x\right\vert ^{2}}{2}u\right)  .
\end{align*}

\begin{proof}
By induction on $k$, we know for any natural number $k$,%
\[
\Delta^{k}\left(  \frac{1+\left\vert x\right\vert ^{2}}{2}u\right)  =k\left(
2k+n-2\right)  \Delta^{k-1}u+2k\sum_{i=1}^{n}x_{i}\Delta^{k-1}\partial
_{i}u+\frac{1+\left\vert x\right\vert ^{2}}{2}\Delta^{k}u.
\]
Then we have%
\begin{align*}
& \Delta\left[  \left(  \frac{1+\left\vert x\right\vert ^{2}}{2}\right)
^{m+1}\Delta^{m}u\right]  +m\left(  m+1\right)  \left(  \frac{1+\left\vert
x\right\vert ^{2}}{2}\right)  ^{m-1}\Delta^{m}u\\
& =\left(  m+1\right)  \left(  2m+n\right)  \left(  \frac{1+\left\vert
x\right\vert ^{2}}{2}\right)  ^{m}\Delta^{m}u+2\left(  m+1\right)  \left(
\frac{1+\left\vert x\right\vert ^{2}}{2}\right)  ^{m}\sum_{i=1}^{n}x_{i}%
\Delta^{m}\partial_{i}u\\
& +\left(  \frac{1+\left\vert x\right\vert ^{2}}{2}\right)  ^{m+1}\Delta
^{m+1}u\\
& =\left(  \frac{1+\left\vert x\right\vert ^{2}}{2}\right)  ^{m}\Delta
^{m+1}\left(  \frac{1+\left\vert x\right\vert ^{2}}{2}u\right)  .
\end{align*}

\end{proof}
\end{lemma}

\begin{lemma}
\label{lemPaneitzexpression}Let $m$ be a natural number and the $2m$th order
operator $A_{2m}$ be given by%
\[
A_{2m}=%
{\displaystyle\prod\limits_{i=0}^{m-1}}
\left[  -\Delta_{S^{n}}+\frac{n\left(  n-2\right)  }{4}-i\left(  i+1\right)
\right]  .
\]
Denote $N$ as the north pole of $S^{n}$ and $\pi_{N}$ as the stereographic
projection from $S^{n}\backslash\left\{  N\right\}  $ to $\mathbb{R}^{n}$,
then for any smooth function $u$ defined on a domain in $\mathbb{R}^{n}$, we
have%
\[
A_{2m}\left(  u\circ\pi_{N}\right)  =\left[  \left(  \frac{2}{1+\left\vert
x\right\vert ^{2}}\right)  ^{-\frac{n+2m}{2}}\left(  -\Delta\right)
^{m}\left(  \left(  \frac{2}{1+\left\vert x\right\vert ^{2}}\right)
^{\frac{n-2m}{2}}u\right)  \right]  \circ\pi_{N}.
\]
In particular, this tells us on $S^{n}$,%
\[
P_{2m}=A_{2m}=%
{\displaystyle\prod\limits_{i=0}^{m-1}}
\left(  -\Delta_{S^{n}}-\left(  i+\frac{n}{2}\right)  \left(  i-\frac{n}%
{2}+1\right)  \right)  .
\]

\begin{proof}
We may identify $S^{n}\backslash\left\{  N\right\}  $ as $\mathbb{R}^{n}$
through $\pi_{N}$, then%
\[
g_{S^{n}}=\frac{4}{\left(  1+\left\vert x\right\vert ^{2}\right)  ^{2}}%
\sum_{i=1}^{n}dx_{i}\otimes dx_{i}.
\]
This implies%
\begin{align*}
\Delta_{S^{n}}u  & =\left(  \frac{2}{1+\left\vert x\right\vert ^{2}}\right)
^{-n}\sum_{i=1}^{n}\partial_{i}\left(  \left(  \frac{2}{1+\left\vert
x\right\vert ^{2}}\right)  ^{n-2}\partial_{i}u\right) \\
& =\left(  \frac{2}{1+\left\vert x\right\vert ^{2}}\right)  ^{-n}\left[
\left(  \frac{2}{1+\left\vert x\right\vert ^{2}}\right)  ^{n-2}\Delta
u-\left(  n-2\right)  \sum_{i=1}^{n}\left(  \frac{2}{1+\left\vert x\right\vert
^{2}}\right)  ^{n-1}x_{i}\partial_{i}u\right] \\
& =\left(  \frac{2}{1+\left\vert x\right\vert ^{2}}\right)  ^{-\frac{n+2}{2}%
}\left[  \left(  \frac{2}{1+\left\vert x\right\vert ^{2}}\right)  ^{\frac
{n-2}{2}}\Delta u-\left(  n-2\right)  \sum_{i=1}^{n}\left(  \frac
{2}{1+\left\vert x\right\vert ^{2}}\right)  ^{\frac{n}{2}}x_{i}\partial
_{i}u\right] \\
& =\left(  \frac{2}{1+\left\vert x\right\vert ^{2}}\right)  ^{-\frac{n+2}{2}%
}\left[  \Delta\left(  \left(  \frac{2}{1+\left\vert x\right\vert ^{2}%
}\right)  ^{\frac{n-2}{2}}u\right)  -\Delta\left(  \frac{2}{1+\left\vert
x\right\vert ^{2}}\right)  ^{\frac{n-2}{2}}\cdot u\right] \\
& =\left(  \frac{2}{1+\left\vert x\right\vert ^{2}}\right)  ^{-\frac{n+2}{2}%
}\Delta\left(  \left(  \frac{2}{1+\left\vert x\right\vert ^{2}}\right)
^{\frac{n-2}{2}}u\right)  +\frac{n\left(  n-2\right)  }{4}u.
\end{align*}
This shows
\[
-\Delta_{S^{n}}u+\frac{n\left(  n-2\right)  }{4}u=\left(  \frac{2}%
{1+\left\vert x\right\vert ^{2}}\right)  ^{-\frac{n+2}{2}}\left(
-\Delta\right)  \left(  \left(  \frac{2}{1+\left\vert x\right\vert ^{2}%
}\right)  ^{\frac{n-2}{2}}u\right)
\]
and verifies the lemma for $m=1$.

Assume the conclusion is true for $m$, then we have%

\begin{align*}
& A_{2\left(  m+1\right)  }u\\
& =\left(  -\Delta_{S^{n}}+\frac{n\left(  n-2\right)  }{4}-m\left(
m+1\right)  \right)  A_{2m}u\\
& =\left(  -\Delta_{S^{n}}+\frac{n\left(  n-2\right)  }{4}\right)  \left[
\left(  \frac{2}{1+\left\vert x\right\vert ^{2}}\right)  ^{-\frac{n+2m}{2}%
}\left(  -\Delta\right)  ^{m}\left(  \left(  \frac{2}{1+\left\vert
x\right\vert ^{2}}\right)  ^{\frac{n-2m}{2}}u\right)  \right] \\
& -m\left(  m+1\right)  \left(  \frac{2}{1+\left\vert x\right\vert ^{2}%
}\right)  ^{-\frac{n+2m}{2}}\left(  -\Delta\right)  ^{m}\left(  \left(
\frac{2}{1+\left\vert x\right\vert ^{2}}\right)  ^{\frac{n-2m}{2}}u\right) \\
& =\left(  \frac{2}{1+\left\vert x\right\vert ^{2}}\right)  ^{-\frac{n+2}{2}%
}\left(  -\Delta\right)  \left[  \left(  \frac{1+\left\vert x\right\vert ^{2}%
}{2}\right)  ^{m+1}\left(  -\Delta\right)  ^{m}\left(  \left(  \frac
{2}{1+\left\vert x\right\vert ^{2}}\right)  ^{\frac{n-2m}{2}}u\right)  \right]
\\
& -m\left(  m+1\right)  \left(  \frac{2}{1+\left\vert x\right\vert ^{2}%
}\right)  ^{-\frac{n+2m}{2}}\left(  -\Delta\right)  ^{m}\left(  \left(
\frac{2}{1+\left\vert x\right\vert ^{2}}\right)  ^{\frac{n-2m}{2}}u\right) \\
& =\left(  \frac{2}{1+\left\vert x\right\vert ^{2}}\right)  ^{-\frac{n+2}{2}%
}\left(  \frac{2}{1+\left\vert x\right\vert ^{2}}\right)  ^{-m}\left(
-\Delta\right)  ^{m+1}\left(  \left(  \frac{2}{1+\left\vert x\right\vert ^{2}%
}\right)  ^{\frac{n-2\left(  m+1\right)  }{2}}u\right) \\
& =\left(  \frac{2}{1+\left\vert x\right\vert ^{2}}\right)  ^{-\frac
{n+2\left(  m+1\right)  }{2}}\left(  -\Delta\right)  ^{m+1}\left(  \left(
\frac{2}{1+\left\vert x\right\vert ^{2}}\right)  ^{\frac{n-2\left(
m+1\right)  }{2}}u\right)  .
\end{align*}
We have used the Lemma \ref{lemidentity1} in the fourth step.
\end{proof}
\end{lemma}

The following basic fact about the Kelvin transformation is an easy corollary
of the above calculations.

\begin{corollary}
\label{corLawOnRn}Assume $u$ is a smooth function. For any Mobius
transformation $\phi$ on $\mathbb{R}^{n}\cup\left\{  \infty\right\}  $,
denote
\[
u_{\phi}=J_{\phi}{}^{\frac{n-2m}{2n}}\cdot u\circ\phi,
\]
here $J_{\phi}$ is the Jacobian of $\phi$, then%
\[
\left(  -\Delta\right)  ^{m}u_{\phi}=J_{\phi}{}^{\frac{n+2m}{2n}}\left(
\left(  -\Delta\right)  ^{m}u\right)  \circ\phi.
\]

\begin{proof}
Since the Mobius transformation group is generated by orthogonal
transformation, translation, dilation and inversion, we only need to verify
the corollary for these special ones. The only nontrivial case is the
inversion. Let $\phi\left(  x\right)  =\frac{x}{\left\vert x\right\vert ^{2}}$
be the inversion map. We may identify $S^{n}\backslash\left\{  N\right\}  $
with $\mathbb{R}^{n}$ through the stereographic projection $\pi_{N}$, then
$\phi^{\ast}g_{S^{n}}=g_{S^{n}}$. It follows that for any function $v$ smooth
away from $0$, we have%
\[
\left(  P_{2m}v\right)  \circ\phi=P_{2m}\left(  v\circ\phi\right)  .
\]
Let
\[
v\left(  x\right)  =\left(  \frac{2}{1+\left\vert x\right\vert ^{2}}\right)
^{-\frac{n-2m}{2}}u\left(  x\right)  ,
\]
then%
\[
\left(  v\circ\phi\right)  \left(  x\right)  =\left(  \frac{2}{1+\left\vert
x\right\vert ^{2}}\right)  ^{-\frac{n-2m}{2}}u_{\phi}\left(  x\right)  .
\]
It follows from Lemma \ref{lemPaneitzexpression} that%
\[
\left(  \left(  P_{2m}v\right)  \circ\phi\right)  \left(  x\right)  =\left(
\frac{1}{\left\vert x\right\vert }\right)  ^{n+2m}\left(  \frac{2}%
{1+\left\vert x\right\vert ^{2}}\right)  ^{-\frac{n+2m}{2}}\left(
-\Delta\right)  ^{m}u\left(  \frac{x}{\left\vert x\right\vert ^{2}}\right)  ,
\]
and%
\[
P_{2m}\left(  v\circ\phi\right)  \left(  x\right)  =\left(  \frac
{2}{1+\left\vert x\right\vert ^{2}}\right)  ^{-\frac{n+2m}{2}}\left(
-\Delta\right)  ^{m}u_{\phi}\left(  x\right)  .
\]
The corollary follows from these two equalities.
\end{proof}
\end{corollary}

It follows from Lemma \ref{lemPaneitzexpression} and Corollary
\ref{corLawOnRn} that

\begin{corollary}
\label{corLawOnSn}Let $u$ be a smooth function on $S^{n}$, $\phi$ be a Mobius
transformation on $S^{n}$,%
\[
u_{\phi}=J_{\phi}^{\frac{n-2m}{2n}}\cdot u\circ\phi,
\]
then%
\[
P_{2m}u_{\phi}=J_{\phi}^{\frac{n+2m}{2n}}\cdot\left(  P_{2m}u\right)
\circ\phi,
\]
and%
\[
E_{2m}\left(  u_{\phi}\right)  =\int_{S^{n}}P_{2m}u_{\phi}\cdot u_{\phi}%
d\mu_{S^{n}}=\int_{S^{n}}P_{2m}u\cdot ud\mu_{S^{n}}=E_{2m}\left(  u\right)  .
\]

\end{corollary}

By Lemma \ref{lemPaneitzexpression}, we may deduce that when $n$ is odd, the
Green's function of $P_{2m}$ at $\xi\in S^{n}$ is equal to%
\begin{equation}
G_{\xi}=\frac{2^{m-n-1}}{\left(  m-1\right)  !%
{\displaystyle\prod\limits_{i=0}^{m}}
\left(  n-2i\right)  \cdot\omega_{n}}\frac{1}{\left(  1+\left\vert \pi_{\xi
}\right\vert ^{2}\right)  ^{\frac{2m-n}{2}}},\label{eqGreenFunction}%
\end{equation}
here $\omega_{n}$ is the volume of the unit ball in $\mathbb{R}^{n}$.

\subsection{Approximation of a Sobolev function by functions vanishing near a
point}

To fully take advantage of the conformal covariant property of the operator
$P_{2m}$, we need to open the punctured sphere as the Euclidean space. Hence
it is useful to understand when a Sobolev function may be approximated by a
sequence of Sobolev functions which vanish near a point.

Let $u$ be a function defined on an open subset of $\mathbb{R}^{n}$, for any
$k>0$, we denote%
\[
D^{k}u=\left(  \partial_{i_{1}i_{2}\cdots i_{k}}u\right)  _{1\leq i_{1}%
,\cdots,i_{k}\leq n}.
\]
We also use the convention $D^{0}u=u$.

\begin{lemma}
\label{lemSmoothing}Assume $1<p<\infty$, $u\in W^{m,p}\left(  B_{1}%
^{n}\right)  $. Let $k$ be the smallest nonnegative integer with $k\geq
m-\frac{n}{p}$. If $D^{j}u\left(  0\right)  =0$ for $0\leq j<k$ (note that the
condition makes sense by the Sobolev embedding theorem, also the condition is
void when $k=0$), then we may find a sequence of smooth functions $u_{i}\in
C^{\infty}\left(  \overline{B}_{1}\right)  $ such that $D^{j}u_{i}\left(
0\right)  =0$ for $0\leq j<k$ and $u_{i}\rightarrow u$ in $W^{m,p}\left(
B_{1}\right)  $.

\begin{proof}
If $k=0$, the conclusion is trivial. Assume $k\geq1$, then it follows from
Sobolev embedding theorem that $W^{m,p}\left(  B_{1}\right)  \subset
C^{k-1}\left(  \overline{B}_{1}\right)  $. First we may find a sequence
$v_{i}\in C^{\infty}\left(  \overline{B}_{1}\right)  $ such that
$v_{i}\rightarrow u$ in $W^{m,p}\left(  B_{1}\right)  $, then the sequence%
\[
u_{i}\left(  x\right)  =v_{i}\left(  x\right)  -\sum_{\left\vert
\alpha\right\vert <k}\frac{\partial^{\alpha}v_{i}\left(  0\right)  }{\alpha
!}x^{\alpha}%
\]
satisfies the requirement in the lemma.
\end{proof}
\end{lemma}

We have the following approximation result, which is a generalization of lemma
2.2 in \cite{HY}.

\begin{proposition}
\label{propVanishNear0}Assume $1<p<\infty$, $u\in W^{m,p}\left(  B_{1}%
^{n}\right)  $. Let $k$ be the smallest nonnegative integer such that $k\geq
m-\frac{n}{p}$. If $D^{j}u\left(  0\right)  =0$ for $0\leq j<k$, then we may
find a sequence $u_{i}\in W^{m,p}\left(  B_{1}^{n}\right)  $ such that
$u_{i}=0$ near the origin, $u_{i}=u$ on $B_{1}\backslash B_{1/2}$ and
$u_{i}\rightarrow u$ in $W^{m,p}\left(  B_{1}^{n}\right)  $.

\begin{proof}
Fix a $\eta\in C^{\infty}\left(  \mathbb{R}^{n}\right)  $ such that $0\leq
\eta\leq1$, $\eta\left(  x\right)  =1$ for $x\in B_{1}$ and $\eta\left(
x\right)  =0$ for $x\in\mathbb{R}^{n}\backslash B_{2}$. For $\lambda>0$, we
let $\eta_{\lambda}\left(  x\right)  =\eta\left(  \frac{x}{\lambda}\right)  $.

First we claim that if $v\in C^{\infty}\left(  \overline{B_{1}}\right)  $ such
that $D^{j}v\left(  0\right)  =0$ for $0\leq j<k$, then we may find a sequence
$v_{i}\in C^{\infty}\left(  \overline{B}_{1}\right)  $ such that
$v_{i}\rightarrow v$ in $W^{m,p}\left(  B_{1}\right)  $ and $v_{i}$ is zero
near the origin. Indeed, let $w_{\varepsilon}=\eta_{\varepsilon}\cdot v$. In
the case when $\left(  m-k\right)  p<n$, we have%
\[
\left\vert D^{m}w_{\varepsilon}\left(  x\right)  \right\vert \leq c\left(
m,n\right)  \sum_{j=0}^{m}\varepsilon^{j-m}\left\vert D^{j}v\left(  x\right)
\right\vert \leq c\left(  m,n,v\right)  \varepsilon^{k-m},
\]
and this implies%
\[
\left\vert D^{m}w_{\varepsilon}\right\vert _{L^{p}\left(  B_{1}\right)  }\leq
c\left(  m,n,v\right)  \varepsilon^{k-m+\frac{n}{p}}\rightarrow0
\]
as $\varepsilon\rightarrow0^{+}$. Hence $v_{\varepsilon}=v-w_{\varepsilon}$ is
the needed approximation function. When $\left(  m-k\right)  p=n$, we only
know $v_{\varepsilon}=v-w_{\varepsilon}$ is bounded in $W^{m,p}\left(
B_{1}\right)  $ and converges to $v$ in $L^{p}\left(  B_{1}\right)  $. Since
$1<p<\infty$, we may find a subsequence $v_{\varepsilon_{i}}\rightharpoonup v$
in $W^{m,p}\left(  B_{1}\right)  $. The claim follows from the standard result
in functional analysis. Indeed, let%
\[
\mathcal{A}=\overline{\operatorname*{co}\left\{  v_{\varepsilon}%
:0<\varepsilon<\frac{1}{16}\right\}  },
\]
here "$\operatorname*{co}$" means the convex hull and the closure is taken in
$W^{m,p}\left(  B_{1}\right)  $. By theorem 2 of chapter 12 in \cite{L}, we
know $\mathcal{A}$ is weakly closed, in particular, $v\in\mathcal{A}$. This
verifies the claim in the case $\left(  m-k\right)  p=n$. We remark that one
may have a constructive proof for this case too.

For any $\varepsilon>0$, by Lemma \ref{lemSmoothing} we may find a $v\in
C^{\infty}\left(  \overline{B_{1}}\right)  $ such that $\left\vert
u-v\right\vert _{W^{m,p}\left(  B_{1}\right)  }\leq\varepsilon$ and
$D^{j}v\left(  0\right)  =0$ for $0\leq j<k$. Then by the above claim we may
find a $\widetilde{v}\in C^{\infty}\left(  \overline{B_{1}}\right)  $ such
that $\widetilde{v}=0$ near the origin and $\left\vert \widetilde
{v}-v\right\vert _{W^{m,p}\left(  B_{1}\right)  }\leq\varepsilon$. Let
$\widetilde{u}=\left(  1-\eta_{1/8}\right)  u+\eta_{1/8}\widetilde{v}$, then
$\widetilde{u}=0$ near the origin, $\widetilde{u}=u$ on $B_{1}\backslash
B_{1/2}$ and $\left\vert \widetilde{u}-u\right\vert _{W^{m,p}\left(
B_{1}\right)  }\leq c\left(  m,p,n\right)  \varepsilon$. The proposition follows.
\end{proof}
\end{proposition}

The same argument will give us the following

\begin{proposition}
Let $u\in W^{m,1}\left(  B_{1}^{n}\right)  $ such that $D^{j}u\left(
0\right)  =0$ for $0\leq j\leq m-n$, then we may find a sequence $u_{i}\in
W^{m,1}\left(  B_{1}^{n}\right)  $ such that $u_{i}=0$ near the origin,
$u_{i}=u$ on $B_{1}\backslash B_{1/2}$ and $u_{i}\rightarrow u$ in
$W^{m,1}\left(  B_{1}^{n}\right)  $.
\end{proposition}

\section{The proof of Theorem \ref{thmmain}\label{sec3}}

\subsection{$n$ is odd and $m=\frac{n+1}{2}$}

In this subsection, we will prove part (i) of Theorem \ref{thmmain}. A crucial
ingredient is the following observation, which should be compared with lemma
7.1 in \cite{HY}.

\begin{lemma}
\label{lemopen1}Assume $n$ is odd and $u\in H^{\frac{n+1}{2}}\left(
S^{n}\right)  $ such that $u\left(  N\right)  =0$, here $N$ is the north pole
of $S^{n}$, then we know
\[
D^{\frac{n+1}{2}}\left(  \sqrt{\frac{1+\left\vert x\right\vert ^{2}}{2}}\cdot
u\left(  \pi_{N}^{-1}\left(  x\right)  \right)  \right)  \in L^{2}\left(
\mathbb{R}^{n}\right)
\]
and%
\[
E_{n+1}\left(  u\right)  =\int_{\mathbb{R}^{n}}\left\vert D^{\frac{n+1}{2}%
}\left(  \sqrt{\frac{1+\left\vert x\right\vert ^{2}}{2}}\cdot u\left(  \pi
_{N}^{-1}\left(  x\right)  \right)  \right)  \right\vert ^{2}dx.
\]
Here $D^{k}f\left(  x\right)  =\left(  \partial_{i_{1}\cdots i_{k}}f\left(
x\right)  \right)  _{1\leq i_{1},\cdots,i_{k}\leq n}$.

\begin{proof}
By Proposition \ref{propVanishNear0} we may find a sequence $u_{i}\in
C^{\infty}\left(  S^{n}\right)  $ such that $u_{i}=0$ near $N$ and
$u_{i}\rightarrow u$ in $H^{\frac{n+1}{2}}\left(  S^{n}\right)  $. By Lemma
\ref{lemPaneitzexpression} we see%
\begin{align*}
& \int_{\mathbb{R}^{n}}\left\vert D^{\frac{n+1}{2}}\left(  \sqrt
{\frac{1+\left\vert x\right\vert ^{2}}{2}}\cdot\left(  u_{i}-u_{j}\right)
\left(  \pi_{N}^{-1}\left(  x\right)  \right)  \right)  \right\vert ^{2}dx\\
& =\int_{\mathbb{R}^{n}}\sqrt{\frac{1+\left\vert x\right\vert ^{2}}{2}}%
\cdot\left(  u_{i}-u_{j}\right)  \left(  \pi_{N}^{-1}\left(  x\right)
\right)  \cdot\left(  -\Delta\right)  ^{\frac{n+1}{2}}\left(  \sqrt
{\frac{1+\left\vert x\right\vert ^{2}}{2}}\cdot\left(  u_{i}-u_{j}\right)
\left(  \pi_{N}^{-1}\left(  x\right)  \right)  \right)  dx\\
& =E_{n+1}\left(  u_{i}-u_{j}\right)  \rightarrow0
\end{align*}
as $i,j\rightarrow\infty$. Hence we may find a vector valued function $F\in
L^{2}\left(  \mathbb{R}^{n}\right)  $ such that%
\[
D^{\frac{n+1}{2}}\left(  \sqrt{\frac{1+\left\vert x\right\vert ^{2}}{2}}\cdot
u_{i}\left(  \pi_{N}^{-1}\left(  x\right)  \right)  \right)  \rightarrow
F\text{ in }L^{2}\left(  \mathbb{R}^{n}\right)  .
\]
This clearly implies%
\[
D^{\frac{n+1}{2}}\left(  \sqrt{\frac{1+\left\vert x\right\vert ^{2}}{2}}\cdot
u\left(  \pi_{N}^{-1}\left(  x\right)  \right)  \right)  =F\in L^{2}\left(
\mathbb{R}^{n}\right)  .
\]
On the other hand, since%
\[
\int_{\mathbb{R}^{n}}\left\vert D^{\frac{n+1}{2}}\left(  \sqrt{\frac
{1+\left\vert x\right\vert ^{2}}{2}}\cdot u_{i}\left(  \pi_{N}^{-1}\left(
x\right)  \right)  \right)  \right\vert ^{2}dx=E_{n+1}\left(  u_{i}\right)  ,
\]
letting $i\rightarrow\infty$, we get
\[
\int_{\mathbb{R}^{n}}\left\vert D^{\frac{n+1}{2}}\left(  \sqrt{\frac
{1+\left\vert x\right\vert ^{2}}{2}}\cdot u\left(  \pi_{N}^{-1}\left(
x\right)  \right)  \right)  \right\vert ^{2}dx=E_{n+1}\left(  u\right)  .
\]

\end{proof}
\end{lemma}

\begin{corollary}
\label{coropen1}Assume $n$ is odd, $u\in H^{\frac{n+1}{2}}\left(
S^{n}\right)  $ and $\xi\in S^{n}$ such that $u\left(  \xi\right)  =0$, then
$E_{n+1}\left(  u\right)  \geq0$. Moreover, $E_{n+1}\left(  u\right)  =0$ if
and only if $u=\operatorname*{const}\cdot\left(  1+\left\vert \pi_{\xi
}\right\vert ^{2}\right)  ^{-1/2}$, here $\pi_{\xi}$ is the stereographic
projection defined at the beginning of Section \ref{sec2}.

\begin{proof}
Without losing of generality, we may assume $\xi=N$. If $E\left(  u\right)  =0
$, then it follows from Lemma \ref{lemopen1} that $\sqrt{\frac{1+\left\vert
x\right\vert ^{2}}{2}}\cdot u\left(  \pi_{N}^{-1}\left(  x\right)  \right)  $
must be a polynomial. On the other hand, since $H^{\frac{n+1}{2}}\left(
S^{n}\right)  \subset C^{\frac{1}{2}}\left(  S^{n}\right)  $ and $u\left(
N\right)  =0$, we see%
\[
\left\vert \sqrt{\frac{1+\left\vert x\right\vert ^{2}}{2}}\cdot u\left(
\pi_{N}^{-1}\left(  x\right)  \right)  \right\vert \leq c\left(  u\right)
\sqrt{\left\vert x\right\vert }%
\]
for $\left\vert x\right\vert $ large. This shows $\sqrt{\frac{1+\left\vert
x\right\vert ^{2}}{2}}\cdot u\left(  \pi_{N}^{-1}\left(  x\right)  \right)
\equiv\operatorname*{const}$. The corollary follows.
\end{proof}
\end{corollary}

Now we are prepared to prove the part (i) of Theorem \ref{thmmain}. The
arguments should be compared to the proof in section 7 of \cite{HY} for
theorem 1.2 there.

\begin{proof}
[Proof of part (i) of Theorem \ref{thmmain}]The key point is to show the
minimizing value of $I_{n+1}$ over $V_{\frac{n+1}{2}}$, $B_{n+1}\left(
S^{n}\right)  $, is reached by some functions. Note that $B_{n+1}\left(
S^{n}\right)  \leq I_{n+1}\left(  1\right)  <0$. Choose a minimizing sequence
$u_{i}\in V_{\frac{n+1}{2}}$ for $I_{n+1}$. By scaling and rotation, we may
assume
\[
\max_{S^{n}}u_{i}=1,\quad\min_{S^{n}}u_{i}=u_{i}\left(  N\right)  .
\]
Here $S$ is the south pole of $S^{n}$. For $i$ large enough, we know
$E_{n+1}\left(  u_{i}\right)  <0$. By the interpolation inequality, we see%
\[
E_{n+1}\left(  u_{i}\right)  \geq c\left\vert u_{i}\right\vert _{H^{\frac
{n+1}{2}}\left(  S^{n}\right)  }^{2}-c\left\vert u_{i}\right\vert
_{L^{2}\left(  S^{n}\right)  }^{2},
\]
this gives us $\left\vert u_{i}\right\vert _{H^{\frac{n+1}{2}}\left(
S^{n}\right)  }\leq c$. After passing to a subsequence, we may find a $u\in
H^{\frac{n+1}{2}}\left(  S^{n}\right)  $ such that $u_{i}\rightharpoonup u$ in
$H^{\frac{n+1}{2}}\left(  S^{n}\right)  $. This implies $u_{i}\rightarrow u$
uniformly on $S^{n}$ and hence $u\geq0$,%
\begin{equation}
\max_{S^{n}}u=1,\quad\min_{S^{n}}u=u\left(  N\right)  .\label{eqNormalization}%
\end{equation}

If $u>0$, then $u_{i}^{-1}\rightarrow u^{-1}$ uniformly on $S^{n}$. Hence
$\left\vert u_{i}^{-1}\right\vert _{L^{2n}\left(  S^{n}\right)  }%
\rightarrow\left\vert u^{-1}\right\vert _{L^{2n}\left(  S^{n}\right)  }$. By
the lower semicontinuity we see%
\[
E_{n+1}\left(  u\right)  \leq\lim\inf_{i\rightarrow\infty}E_{n+1}\left(
u_{i}\right)  ,
\]
hence we see%
\[
B_{n+1}\left(  S^{n}\right)  \leq\left\vert u^{-1}\right\vert _{L^{2n}\left(
S^{n}\right)  }^{2}E_{n+1}\left(  u\right)  \leq\lim\inf_{i\rightarrow\infty
}\left\vert u^{-1}\right\vert _{L^{2n}\left(  S^{n}\right)  }^{2}%
E_{n+1}\left(  u_{i}\right)  =B_{n+1}\left(  S^{n}\right)  ,
\]
and $u$ is a minimizer.

If $u$ vanishes at some point, say $\xi\in S^{n}$. It follows from lower
semicontinuity that $E_{n+1}\left(  u\right)  \leq0$. By Corollary
\ref{coropen1} we see $u=c\cdot\left(  1+\left\vert \pi_{\xi}\right\vert
^{2}\right)  ^{-1/2}$. Using (\ref{eqNormalization}), we see $c=1,\xi=N$. In
particular, $u_{i}\left(  S\right)  \rightarrow1$ as $i\rightarrow\infty$.
Denote%
\[
\lambda_{i}=\frac{u_{i}\left(  S\right)  }{u_{i}\left(  N\right)  }%
\rightarrow\infty.
\]
Using the notations in Corollary \ref{corLawOnSn} and the beginning of Section
\ref{sec2}, we let%
\[
v_{i}=\left(  u_{i}\right)  _{\sigma_{N,\lambda_{i}}}=\left(  \frac
{1+\lambda_{i}^{2}\left\vert \pi_{N}\right\vert ^{2}}{\lambda_{i}\left(
1+\left\vert \pi_{N}\right\vert ^{2}\right)  }\right)  ^{1/2}\cdot u_{i}%
\circ\sigma_{N,\lambda_{i}}.
\]
Then $v_{i}$ is still a minimizing sequence for $I_{n+1}$ with $v_{i}\left(
N\right)  =v_{i}\left(  S\right)  $. Let $\nu_{i}=\max_{S^{n}}v_{i}$,
$w_{i}=\frac{v_{i}}{\nu_{i}}$, then $w_{i}$ is a minimizing sequence and after
passing to a subsequence, we may find a $w\in H^{\frac{n+1}{2}}\left(
S^{n}\right)  $ such that $w_{i}\rightharpoonup w$ in $H^{\frac{n+1}{2}%
}\left(  S^{n}\right)  $. We claim $w>0$. Indeed, if this is not the case,
then for some $\xi\in S^{n}$, $w\left(  \xi\right)  =0$. Argue as before we
see $w=\left(  1+\left\vert \pi_{\xi}\right\vert ^{2}\right)  ^{-1/2}$. Since
$w\left(  S\right)  =w\left(  N\right)  $, we see $\xi\neq N,S$. In
particular,%
\[
w_{i}\left(  S\right)  =\frac{u_{i}\left(  S\right)  }{\nu_{i}\sqrt
{\lambda_{i}}}\rightarrow w\left(  S\right)  >0.
\]
On $S^{n}\backslash\left\{  S,N\right\}  $, we have%
\begin{align*}
w_{i}  & \geq\frac{u_{i}\left(  N\right)  }{\nu_{i}}\left(  \frac
{1+\lambda_{i}^{2}\left\vert \pi_{N}\right\vert ^{2}}{\lambda_{i}\left(
1+\left\vert \pi_{N}\right\vert ^{2}\right)  }\right)  ^{1/2}\\
& =\frac{u_{i}\left(  S\right)  }{\nu_{i}\sqrt{\lambda_{i}}}\left(
\frac{\lambda_{i}^{-2}+\left\vert \pi_{N}\right\vert ^{2}}{1+\left\vert
\pi_{N}\right\vert ^{2}}\right)  ^{1/2}\rightarrow\frac{w\left(  S\right)
\left\vert \pi_{N}\right\vert }{\left(  1+\left\vert \pi_{N}\right\vert
^{2}\right)  ^{1/2}},
\end{align*}
this implies $w>0$ on $S^{n}$ and contradicts with our assumption. Hence $w>0
$ and it is a minimizer.

Assume $u$ is a minimizer for $I_{n+1}$ in $V_{\frac{n+1}{2}}$, then for some
positive constant $c$, we have%
\[
u\in C^{\infty}\left(  S^{n}\right)  ,\text{ }u>0\text{ and }P_{n+1}%
u=-cu^{-2n-1}\text{ on }S^{n}.
\]
Using the Green's function of $P_{2m}$ written down at the end of Section
\ref{secP2mExpression}, we see for some $c>0$,%
\[
u\left(  \xi\right)  =c\int_{S^{n}}\frac{u\left(  \zeta\right)  ^{-2n-1}%
}{\left(  1+\left\vert \pi_{\xi}\left(  \zeta\right)  \right\vert ^{2}\right)
^{1/2}}d\mu_{S^{n}}\left(  \zeta\right)  \quad\text{for any }\xi\in S^{n}.
\]
Let
\[
v\left(  x\right)  =\sqrt{\frac{1+\left\vert x\right\vert ^{2}}{2}}\cdot
u\left(  \pi_{N}^{-1}\left(  x\right)  \right)  ,
\]
then%
\[
\left(  \pi_{N}^{-1}\right)  ^{\ast}\left(  u^{-4}g_{S^{n}}\right)
=v^{-4}g_{\mathbb{R}^{n}}.
\]
Moreover, it follows from the integral equation of $u$ that for some $c>0$,%
\[
v\left(  x\right)  =c\int_{\mathbb{R}^{n}}\left\vert x-y\right\vert v\left(
y\right)  ^{-2n-1}dy.
\]
It follows from theorem 1.5 of \cite{Li} (proved by the method of moving
spheres, a variation of the method of moving planes \cite{GNN}, see also
\cite{CLO} for the method of moving planes for integral equations) that%
\[
v\left(  x\right)  =c\left(  \frac{1+\lambda^{2}\left\vert x-x_{0}\right\vert
^{2}}{2\lambda}\right)  ^{1/2}%
\]
for some $c>0,\lambda>0$ and $x_{0}\in\mathbb{R}^{n}$. It follows that for
some Mobius transformation $\phi$ on $S^{n}$, we have $u=cJ_{\phi}^{-\frac
{1}{2n}}$. Using Corollary \ref{corLawOnSn} we see that $I_{n+1}\left(
cJ_{\phi}^{-\frac{1}{2n}}\right)  =I_{n+1}\left(  1\right)  $. Hence $1$ is a
minimizer of $I_{n+1}$ and all the minimizers are of the form $cJ_{\phi
}^{-\frac{1}{2n}}$ for some Mobius transformation $\phi$.
\end{proof}

\subsection{$n$ is odd and $m=\frac{n+3}{2}$}

The argument for the part (ii) of Theorem \ref{thmmain} goes along the similar
line as for part (i). We will only explain when the proof is different. First
we have

\begin{lemma}
\label{lemopen2}Assume $n$ is odd and $u\in H^{\frac{n+3}{2}}\left(
S^{n}\right)  $ such that $u\left(  N\right)  =0$ and $du\left(  N\right)
=0$, here $N$ is the north pole of $S^{n}$, then we know
\[
D^{\frac{n+3}{2}}\left(  \left(  \frac{1+\left\vert x\right\vert ^{2}}%
{2}\right)  ^{3/2}\cdot u\left(  \pi_{N}^{-1}\left(  x\right)  \right)
\right)  \in L^{2}\left(  \mathbb{R}^{n}\right)
\]
and%
\[
E_{n+3}\left(  u\right)  =\int_{\mathbb{R}^{n}}\left\vert D^{\frac{n+3}{2}%
}\left(  \left(  \frac{1+\left\vert x\right\vert ^{2}}{2}\right)  ^{3/2}\cdot
u\left(  \pi_{N}^{-1}\left(  x\right)  \right)  \right)  \right\vert ^{2}dx.
\]

\end{lemma}

Similar to Lemma \ref{lemopen1}, this lemma follows from an approximation
argument using Proposition \ref{propVanishNear0} and Lemma
\ref{lemPaneitzexpression}.

\begin{corollary}
\label{coropen2}Assume $n$ is odd, $u\in H^{\frac{n+3}{2}}\left(
S^{n}\right)  ,$ $u\geq0$ and $\xi\in S^{n}$ such that $u\left(  \xi\right)
=0$, then $E_{n+3}\left(  u\right)  \geq0$. Moreover, $E_{n+3}\left(
u\right)  =0$ if and only if $u=\operatorname*{const}\cdot\left(  1+\left\vert
\pi_{\xi}\right\vert ^{2}\right)  ^{-3/2}$.

\begin{proof}
Without losing of generality, we may assume $\xi=N$. Since $u\in H^{\frac
{n+3}{2}}\left(  S^{n}\right)  \subset C^{1,1/2}\left(  S^{n}\right)  $, we
see $du\left(  N\right)  =0$. It follows from Lemma \ref{lemopen2} that
$E_{n+3}\left(  u\right)  \geq0$. Moreover, if $E_{n+3}\left(  u\right)  =0$,
then
\[
D^{\frac{n+3}{2}}\left(  \left(  \frac{1+\left\vert x\right\vert ^{2}}%
{2}\right)  ^{3/2}\cdot u\left(  \pi_{N}^{-1}\left(  x\right)  \right)
\right)  =0.
\]
This implies $\left(  \frac{1+\left\vert x\right\vert ^{2}}{2}\right)
^{3/2}\cdot u\left(  \pi_{N}^{-1}\left(  x\right)  \right)  $ must be a
polynomial. Since%
\[
\left\vert \left(  \frac{1+\left\vert x\right\vert ^{2}}{2}\right)
^{3/2}\cdot u\left(  \pi_{N}^{-1}\left(  x\right)  \right)  \right\vert \leq
c\left(  u\right)  \left\vert x\right\vert ^{3/2}%
\]
when $\left\vert x\right\vert $ is large, we see%
\[
\left(  \frac{1+\left\vert x\right\vert ^{2}}{2}\right)  ^{3/2}\cdot u\left(
\pi_{N}^{-1}\left(  x\right)  \right)  =c_{0}+\sum_{i=1}^{n}c_{i}x_{i}.
\]
It follows from the fact $u\geq0$ that $c_{0}\geq0$ and $c_{i}=0$ for $1\leq
i\leq n$. Hence%
\[
u=c_{0}\left(  \frac{1+\left\vert x\right\vert ^{2}}{2}\right)  ^{-3/2}.
\]

\end{proof}
\end{corollary}

\begin{proof}
[Sketch of the proof of part (ii) of Theorem \ref{thmmain}]The key point is to
show the minimizing value $B_{n+3}\left(  S^{n}\right)  $ is reached at some
function, then one may use the theorem 1.5 of \cite{Li} (see also closely
related results in \cite{CLO}) and Corollary \ref{corLawOnSn} to conclude that
$u=1$ is a minimizer.

Let $u_{i}$ be a minimizing sequence for $I_{n+3}$ in $V_{\frac{n+3}{2}}$, by
scaling and rotation we may assume
\[
\max_{S^{n}}u_{i}=1\text{ and }\min_{S^{n}}u_{i}=u_{i}\left(  N\right)  .
\]
Then since
\[
\left\vert u_{i}^{-1}\right\vert _{L^{\frac{2n}{3}}\left(  S^{n}\right)  }%
^{2}E_{n+3}\left(  u_{i}\right)  \leq c,
\]
we see $E_{n+3}\left(  u_{i}\right)  \leq c$. By coercivity we see $\left\vert
u_{i}\right\vert _{H^{\frac{n+3}{2}}\left(  S^{n}\right)  }\leq c$. After
passing to a subsequence, we may find a $u\in H^{\frac{n+3}{2}}\left(
S^{n}\right)  $ such that $u_{i}\rightharpoonup u$ in $H^{\frac{n+3}{2}%
}\left(  S^{n}\right)  $. Then $u_{i}\rightarrow u$ uniformly. We have $u\geq0
$ and%
\begin{equation}
\max_{S^{n}}u=1,\quad\min_{S^{n}}u=u\left(  N\right)
.\label{eqNormalization2}%
\end{equation}

If $u>0$, then it is a minimizer as before.

If $u\left(  \xi\right)  =0$ for some $\xi\in S^{n}$, then since $u\geq0$ and
$u\in C^{1,1/2}\left(  S^{n}\right)  $, we see $du\left(  \xi\right)  =0$ and%
\[
\left\vert u\left(  \zeta\right)  \right\vert \leq c\left(  u\right)
d_{S^{n}}\left(  \zeta,\xi\right)  ^{3/2}\text{ for }\zeta\in S^{n}.
\]
Here $d_{S^{n}}\left(  \zeta,\xi\right)  $ is the geodesic distance on $S^{n}
$ with standard metric. This implies $\left\vert u^{-1}\right\vert
_{L^{\frac{2n}{3}}\left(  S^{n}\right)  }=\infty$. It follows from Fatou's
lemma that%
\[
\infty=\left\vert u^{-1}\right\vert _{L^{\frac{2n}{3}}\left(  S^{n}\right)
}\leq\lim\inf_{i\rightarrow\infty}\left\vert u_{i}^{-1}\right\vert
_{L^{\frac{2n}{3}}\left(  S^{n}\right)  },
\]
hence $\left\vert u_{i}^{-1}\right\vert _{L^{\frac{2n}{3}}\left(
S^{n}\right)  }\rightarrow\infty$ as $i\rightarrow\infty$. Using lower
semicontinuity we see $E\left(  u\right)  \leq0$. It follows from Corollary
\ref{coropen2} that%
\[
u=c\left(  1+\left\vert \pi_{\xi}\right\vert ^{2}\right)  ^{-3/2}.
\]
In view of (\ref{eqNormalization2}), we see $c=1$ and $\xi=N$, hence
$u=\left(  1+\left\vert \pi_{N}\right\vert ^{2}\right)  ^{-3/2}$. Now we may
proceed to renormalize the minimizing sequence as in the proof of part (i).
\end{proof}

\subsection{$n$ is odd and $m\geq\frac{n+5}{2}$}

The arguments presented in the previous two subsections do not work well for
the case when $n$ is odd and $m\geq\frac{n+5}{2}$. The main problem is that we
do not get enough "vanishing condition" when the weak limit of minimizing
sequence touches zero. This looks like a technical point. But in fact, it is
essential, we will show no minimizer exists at all when $m$ becomes this larger.

\begin{proof}
[Proof of part (iii) of Theorem \ref{thmmain}]First we observe that it follows
from theorem 1.5 of \cite{Li} (see also closely related results in \cite{CLO})
that any critical point of $I_{2m}$ over $V_{m}$ must be of the form
$cJ_{\phi}^{\frac{n-2m}{2n}}$ for some $c>0 $ and Mobius transformation $\phi
$. In view of Corollary \ref{corLawOnSn}, to show all of them are unstable, we
only need to show $u=1$ is unstable. Calculation shows the second variation of
$I_{2m}$ at $u=1$, namely $H$, is given by%
\begin{align*}
& \frac{1}{2}\left(  \mu_{S^{n}}\left(  S^{n}\right)  \right)  ^{-\frac
{2m-n}{n}}H\left(  \varphi\right) \\
& =E_{2m}\left(  \varphi\right)  +\frac{2m+n}{2m-n}P_{2m}1\cdot\int_{S^{n}%
}\varphi^{2}d\mu_{S^{n}}-\frac{4m}{2m-n}\frac{P_{2m}1}{\mu_{S^{n}}\left(
S^{n}\right)  }\left(  \int_{S^{n}}\varphi d\mu_{S^{n}}\right)  ^{2}%
\end{align*}
for any $\varphi\in H^{m}\left(  S^{n}\right)  $. The corresponding
self-adjoint operator is given by%
\[
\mathcal{A}\varphi=P_{2m}\varphi+\frac{2m+n}{2m-n}P_{2m}1\cdot\varphi
-\frac{4m}{2m-n}\frac{P_{2m}1}{\mu_{S^{n}}\left(  S^{n}\right)  }\int_{S^{n}%
}\varphi d\mu_{S^{n}}.
\]

When $m-\frac{n+5}{2}$ is even, let $h_{2}$ be any harmonic homogeneous
polynomial of degree $2$, then%
\[
\mathcal{A}h_{2}=2m%
{\displaystyle\prod\limits_{i=0}^{m}}
\left(  \frac{n}{2}+i\right)
{\displaystyle\prod\limits_{i=1}^{m-2}}
\left(  \frac{n}{2}-i\right)  \cdot h_{2}.
\]
It gives us a negative eigenvalue.

When $m-\frac{n+5}{2}$ is odd, let $h_{3}$ be any harmonic homogeneous
polynomial of degree $3$, then%
\begin{align*}
& \mathcal{A}h_{3}\\
& =\left[  \left(  m+\frac{n}{2}+1\right)  \left(  m+\frac{n}{2}+2\right)
-\left(  m-\frac{n}{2}-2\right)  \left(  m-\frac{n}{2}-1\right)  \right] \\
& \cdot%
{\displaystyle\prod\limits_{i=0}^{m}}
\left(  \frac{n}{2}+i\right)
{\displaystyle\prod\limits_{i=1}^{m-3}}
\left(  \frac{n}{2}-i\right)  \cdot h_{3}.
\end{align*}

Again, it gives us a negative eigenvalue.
\end{proof}

\subsection{$n$ is even}

The case when the dimension is even is very different form the odd dimension.
In fact, in this case, the variational problem (\ref{eqmain}) becomes trivial.

\begin{proof}
[Proof of part (iv) of Theorem \ref{thmmain} ]This follows from the formula of
$P_{2m}$ (cf. Lemma \ref{lemPaneitzexpression}) and the fact that the
eigenvalues of $-\Delta_{S^{n}}$ are given by $\alpha\left(  \alpha
+n-1\right)  $, $\alpha\in\mathbb{Z}_{+}$, with corresponding eigenfunctions
given by harmonic homogeneous polynomials of degree $\alpha$.
\end{proof}

\section{Another approach to Theorem \ref{thmmain}\label{sec4}}

In deriving an upper bound for the eigenvalue of an arbitrary metric on
$S^{2}$, Hersch used the conformal invariance property of the Dirichlet energy
to choose suitable test functions through a trick which became popular later
and is known as the barycenter analysis (see p142 of \cite{SY}). Such kind of
trick was used in \cite{ACW} for the proof of (\ref{eqP2S1}) and more recently
in \cite{NZ} for (\ref{eqP4S1}). In this section, we will combine this trick
with Lemma \ref{lemopen1} and \ref{lemopen2} to give another approach for
Theorem \ \ref{thmmain}.

For any $a\in B_{1}^{n+1}$, we have a smooth diffeomorphism from $\overline
{B}_{1}^{n+1}$ to itself given by%
\[
\sigma_{a}\left(  z\right)  =\frac{\left(  1-\left\vert a\right\vert
^{2}\right)  z-\left(  \left\vert z\right\vert ^{2}-2a\cdot z+1\right)
a}{\left\vert a\right\vert ^{2}\left\vert z\right\vert ^{2}-2a\cdot z+1}\text{
for }z\in\overline{B}_{1}^{n+1}.
\]
Note that $\sigma_{a}\left(  a\right)  =0$, $\sigma_{a}^{-1}=\sigma_{-a}$ and
for $a\neq0$,%
\[
\left.  \sigma_{a}\right\vert _{S^{n}}=\sigma_{\frac{a}{\left\vert
a\right\vert },\frac{1-\left\vert a\right\vert }{1+\left\vert a\right\vert }}.
\]
Let $2m>n$, $u\in C\left(  S^{n},\mathbb{R}\right)  $ be a strictly positive
function. For $a\in B_{1}^{n+1}$, let $u_{\sigma_{a}}=J_{\sigma_{a}}%
^{\frac{n-2m}{2n}}\cdot u\circ\sigma_{a}$ and $C\left(  a\right)  =\int
_{S^{n}}u_{\sigma_{a}}\left(  \zeta\right)  \zeta d\mu_{S^{n}}\left(
\zeta\right)  $, then%
\[
C\left(  a\right)  =\left(  2\lambda\right)  ^{\frac{n-2m}{2}}\int_{S^{n}%
}u\left(  \sigma_{a}\left(  \zeta\right)  \right)  \left[  1-\frac{\zeta\cdot
a}{\left\vert a\right\vert }+\lambda^{2}\left(  1+\frac{\zeta\cdot
a}{\left\vert a\right\vert }\right)  \right]  ^{\frac{2m-n}{2}}\zeta
d\mu_{S^{n}}\left(  \zeta\right)  ,
\]
here $\lambda=\frac{1-\left\vert a\right\vert }{1+\left\vert a\right\vert }$.
In particular%
\[
\left(  2\lambda\right)  ^{\frac{2m-n}{2}}C\left(  a\right)  \rightarrow
\int_{S^{n}}u\left(  -\xi\right)  \left(  1-\zeta\cdot\xi\right)
^{\frac{2m-n}{2}}\zeta d\mu_{S^{n}}\left(  \zeta\right)  =-c\left(
m,n\right)  u\left(  -\xi\right)  \xi
\]
as $a\rightarrow\xi\in S^{n}$. Since $C$ is continuous on $B_{1}^{n+1}$, it
follows from winding number argument that for some $a\in B_{1}^{n+1}$,
$C\left(  a\right)  =0$.

Now we may sketch a somewhat different argument for part (i) and (ii) of
Theorem \ref{thmmain}. We restrict ourselves to part (i) since the argument
for part (ii) is very similar. Let $u_{i}$ be a minimizing sequence for
$I_{n+1}$ over $V_{\frac{n+1}{2}}$, we may find $a_{i}\in B_{1}^{n+1}$ such
that $\int_{S^{n}}\left(  u_{i}\right)  _{\sigma_{a_{i}}}\left(  \zeta\right)
\zeta d\mu_{S^{n}}\left(  \zeta\right)  =0$. Since $I_{n+1}\left(
u_{i}\right)  =I_{n+1}\left(  \left(  u_{i}\right)  _{\sigma_{a_{i}}}\right)
$, we may assume $\int_{S^{n}}u_{i}\left(  \zeta\right)  \zeta d\mu_{S^{n}%
}\left(  \zeta\right)  =0$. By scaling and rotation we may also assume
$\max_{S^{n}}u_{i}=1$ and $\min_{S^{n}}u_{i}=u_{i}\left(  N\right)  $. The
same argument as in Section \ref{sec3} shows for some $u\in H^{\frac{n+1}{2}%
}\left(  S^{n}\right)  $, we have $u_{i}\rightharpoonup u$ in $H^{\frac
{n+1}{2}}\left(  S^{n}\right)  $. We only need to show $u>0$ on $S^{n}$. If
$u$ touches zero somewhere, then as in Section \ref{sec3}, we see $u=\left(
1+\left\vert \pi_{N}\right\vert ^{2}\right)  ^{-1/2}$. On the other hand, it
follows from $\int_{S^{n}}u_{i}\left(  \zeta\right)  \zeta d\mu_{S^{n}}\left(
\zeta\right)  =0$ that $\int_{S^{n}}u\left(  \zeta\right)  \zeta d\mu_{S^{n}%
}\left(  \zeta\right)  =0$. But $\int_{S^{n}}\left(  1+\left\vert \pi
_{N}\left(  \zeta\right)  \right\vert ^{2}\right)  ^{-1/2}\zeta\neq0$, this
gives us a contradiction. Hence $u$ never touches zero and it must be a
minimizer. The remaining argument is the same as in Section \ref{sec3}.

\section{Further remarks}

Recently in \cite{NZ}, an argument for (\ref{eqP4S1}) is given based on the
observation that $P_{4}$ is positive definite on the $L^{2}$ orthogonal
complement of the restrictions of linear functions on $S^{1}$. Such kind of
argument works in higher dimension for part (ii) of Theorem \ref{thmmain} too.
Indeed, we note that if $n$ is even, then%
\[
P_{n+3}=\left(  -\Delta_{S^{n}}-\frac{n-\frac{1}{2}}{2}\right)  \left(
-\Delta_{S^{n}}-\frac{3\left(  n+\frac{1}{2}\right)  }{2}\right)
{\displaystyle\prod\limits_{i=0}^{\frac{n-3}{2}}}
\left(  -\Delta_{S^{n}}+\left(  i+\frac{n}{2}\right)  \left(  \frac{n-2}%
{2}-i\right)  \right)  .
\]
Using the fact that the eigenvalues of $-\Delta_{S^{n}}$ are given by
$\alpha\left(  \alpha+n-1\right)  $, $\alpha\in\mathbb{Z}_{+}$, with
corresponding eigenfunctions given by harmonic homogeneous polynomials of
degree $\alpha$, we see $P_{n+3}$ is positive definite on the $L^{2}$
orthogonal complement of linear functions. Assume $u_{i}\in V_{\frac{n+3}{2}}
$ is a minimizing sequence of $I_{n+3}$. Without losing of generality, we may
assume $\max_{S^{n}}u_{i}=1$ and $u_{i}$ is perpendicular to linear functions.
By the arguments in Section \ref{sec3}, we may find $u\in H^{\frac{n+3}{2}%
}\left(  S^{n}\right)  $ such that $u_{i}\rightharpoonup u$ in $H^{\frac
{n+3}{2}}\left(  S^{n}\right)  $. If $u$ touches $0$ somewhere, then as in
Section \ref{sec3} we know $E_{n+3}\left(  u\right)  \leq0$ and $\max_{S^{n}%
}u=1$. This contradicts with the fact that $P_{n+3}$ is strictly positive
definite on the orthogonal complement of linear functions. Hence $u$ does not
touch zero and it is a minimizer.

\end{document}